\documentclass[12pt]{article}
\usepackage{amsfonts}
\usepackage{makeidx}
\usepackage{amssymb}
\usepackage{graphicx}
\usepackage{amsmath}
\usepackage[T1]{fontenc}
\usepackage{amsfonts}

\newtheorem{theorem}{Theorem}

\newtheorem{adefinition}[theorem]{Definition}
\newenvironment{definition}{\begin{adefinition}\rm}{\end{adefinition}}
\newtheorem{aexample}[theorem]{Example}

\newtheorem{lemma}[theorem]{Lemma}

\newtheorem{proposition}[theorem]{Proposition}
\newtheorem{aremark}[theorem]{Remark}
\newenvironment{remark}{\begin{aremark}\rm}{\end{aremark}}

\numberwithin{equation}{section} \numberwithin{theorem}{section}
\newenvironment{proof}[1][Proof]{\textbf{#1.} }{\ \rule{0.5em}{0.5em}}
\makeatother

\begin{document}

\title{On the Limiting Empirical Measure of the sum of rank one matrices with log-concave distribution}
\author{A. Pajor$^1$ \\
Department of Mathematics, University Paris-Est, \\
Marne la Vall\'{e}e, France \and L. Pastur$^1$\\
Theoretical Department, Institute for Low Temperatures\\
Kharkiv, Ukraine}
\date{}
\maketitle \footnotetext[1]{Research supported in part by the
European Network PHD, MCRN-511953}

\begin{abstract}
We consider $n\times n$ real symmetric and hermitian random matrices $H_{n,m}
$ equals the sum of a non-random matrix $H_{n}^{(0)}$ matrix and the sum of $%
m$ rank-one matrices determined by $m$ i.i.d. isotropic random
vectors with log-concave probability law and i.i.d. random
amplitudes $\{\tau _{\alpha }\}_{\alpha =1}^{m}$. This is a
generalization of the case of vectors uniformly distributed over the
unit sphere, studied in \cite{Ma-Pa:67}. We prove as in
\cite{Ma-Pa:67} that if $n\rightarrow \infty ,\;m\rightarrow \infty
,\;m/n\rightarrow c\in \lbrack 0,\infty )$ and that the empirical
eigenvalue measure of $H_{n}^{(0)}$ converges weakly, then the
empirical eigenvalue measure of $H_{n,m}$ converges in probability
to a non-random limit, found in  \cite{Ma-Pa:67}.
\end{abstract}

\section{Introduction: Problem and Main Result}

Let $\{Y_{\alpha }\}_{\alpha =1}^{m}$ be i.i.d. random vectors of $\mathbb{R}%
^{n}$ (or $\mathbb{C}^{n}$), and $\{\tau _{\alpha }\}_{\alpha
=1}^{m}$ be
i.i.d. random variables with common probability law $\sigma $. Set%
\begin{equation}
M_{n,m}=\sum_{\alpha =1}^{m}\tau _{\alpha }L_{Y_{\alpha }},  \label{Mnm}
\end{equation}%
where$\;L_{Y}=Y\otimes Y$ is the rank-one matrix, corresponding to $Y\in \mathbb{R}%
^{n}( $or $\mathbb{C}^{n})$ and defined as%
\begin{equation}
L_{Y}X=(X,Y)Y,\;\forall X\in \mathbb{R}^{n}(\mathbb{C}^{n}),  \label{LYX}
\end{equation}%
with $(\;,\;)$ denoting the standard euclidian (or hermitian) scalar product
in $\mathbb{R}^{n}($or $\mathbb{C}^{n})$.

Let also $H_{n}^{(0)}$ be a real symmetric (or hermitian) $n\times n$
matrix. We then consider the real symmetric (or hermitian) $n\times n$ random%
\begin{equation}
H_{n,m}=H_{n}^{(0)}+M_{n,m}.  \label{Hnm}
\end{equation}%
Denote%
\begin{equation}
-\infty <\lambda _{1}\leq ...\leq \lambda _{n}<\infty \label{evHnm}
\end{equation}%
the eigenvalues of $H_{n,m}$ and introduce their Normalized Counting (or
empirical) Measure $N_{n,m}$, setting for any interval $\Delta \subset
\mathbb{R}^{n}$%
\begin{equation}
N_{n,m}(\Delta )=\mathrm{Card}\{l\in \lbrack 1,n]:\lambda _{l}\in \Delta
\}/n.  \label{NCM}
\end{equation}%
Likewise, we define the Normalized Counting Measure $N_{n}^{(0)}$ of
eigenvalues $\{\lambda _{l}^{(0)}\}_{l=1}^{n}$ of $H_{n}^{(0)}$%
\begin{equation}
N_{n}^{(0)}(\Delta )=\mathrm{Card}\{l\in \lbrack 1,n]:\lambda _{l}^{(0)}\in
\Delta \}/n  \label{NCM0}
\end{equation}%
and we assume that the sequence $\{N_{n}^{(0)}\}$ converges weakly to a
probability measure $N^{(0)}$:%
\begin{equation}
\lim_{n\rightarrow \infty }N_{n}^{(0)}=N^{(0)}.  \label{N0N}
\end{equation}%
It was shown in \cite{Ma-Pa:67} that if $\{Y_{\alpha }\}_{\alpha =1}^{m}$
are uniformly distributed over the unit sphere of $\mathbb{R}^{n}(\mathbb{C}%
^{n})$ and
\begin{equation}
n\rightarrow \infty ,\;m\rightarrow \infty ,\;m/n\rightarrow c\in \lbrack
0,\infty ),  \label{nmc}
\end{equation}%
then there exists a non-random probability measure $N$
($N(\mathbb{R})=1$) such that for any interval $\Delta \subset
\mathbb{R}$ we have the convergence in probability
\begin{equation}
\lim_{n\rightarrow \infty ,\;m\rightarrow \infty ,\;m/n\rightarrow
c}N_{n,m}(\Delta )=N(\Delta ).  \label{NnmN}
\end{equation}%
The limiting non-random measure $N$ can be found as follows. Introduce the
Stieltjes transform
\begin{equation}
{f}^{(0)}(z)=\int_{\mathbb{R}}\frac{N^{(0)}(d\lambda )}{\lambda -z},\;\Im
z\neq 0  \label{f0}
\end{equation}%
of the measure $N^{(0)}$ of (\ref{N0N}) and the Stieltjes transform%
\begin{equation}
f(z)=\int_{\mathbb{R}}\frac{N(d\lambda )}{\lambda -z},\;\Im z\neq 0
\label{f}
\end{equation}%
of the measure $N$ of (\ref{NnmN}). Then $f$ is uniquely determined by the
functional equation%
\begin{equation}
f(z)=f^{(0)}\left( z-c\int_{\mathbb{R}}\frac{\tau \sigma (d\tau )}{1+\tau
f(z)}\right) ,  \label{MPE}
\end{equation}%
considered in the class of functions analytic in $\mathbb{C\setminus \mathbb{%
R}}$ and such that $\Im f(z)\Im z\geq 0,\;\Im z\neq 0$. Since the Stieltjes
transform determines $N$ uniquely by the formula%
\begin{equation}
\lim_{\varepsilon \rightarrow 0^{+}}\frac{1}{\pi }\int_{\mathbb{R}}\varphi
(\lambda )\Im f(\lambda +i\varepsilon )d\lambda =\int_{\mathbb{R}}\varphi
(\lambda )N(d\lambda ),  \label{SP}
\end{equation}%
valid for any continuous function of compact support, (\ref{MPE}) determines
$N$ uniquely.

The same result is valid if the components $\{Y_{\alpha j}\}_{j=1}^{n}$ of $%
Y_{\alpha },\;\alpha =1,...,m$ are i.i.d. random variables of zero mean and
of unit variance. A particular case of this for $H_{n}^{(0)}=0,\;\tau
_{\alpha }=1,\;\alpha =1,...,m$ and Gaussian $\{Y_{\alpha j}\}_{j=1}^{n}$ is
known since the 30th in statistics as the Wishart matrix (see e.g. \cite%
{Mui:82}). The same random matrix \ appears also in the local theory of
Banach spaces or so-called asymptotic convex geometry (see e.g. \cite%
{{Da-Sz:01},{Sz:06}}). One particular and important case that enters
this framework is the study of some geometric parameters
associated to i.i.d. random points
$Y_{\alpha },\;\alpha =1,...,m$ uniformly distributed on a convex
body in $\mathbb{R}^n$ and
the asymptotic geometry of the random convex polytope generated by these points
(see e.g. \cite%
{{Bou:96}, {Gi-Mi:00}, {Gi-Ts:03},{Li-Ru-Paj-Tom:05}, {Au:07}}).

In this paper we prove the same result for the case where the common
probability law of the i.i.d. vectors $Y_{\alpha },\;\alpha =1,...,m$ is
isotropic and log-concave. A preliminary non-published result was obtained
in 2004 by the authors when the vectors are random points uniformly
distributed in the unit ball $\{\sum_1^n |x_i|^p\le 1\}$. The latter case
was also obtained by a different approach in \cite{Au:06}.  Here are the
corresponding definitions.

\begin{definition}
[Isotropic vectors] A random real vector $Y\in \mathbb{R}^{n}$ is called
isotropic if
\begin{equation}
\mathbf{E}\{(Y,X)\}=0\quad \mathrm{and}\quad \mathbf{E}\{|(Y,X)|^{2}%
\}=n^{-1}|X|^{2},\;\forall X\in \mathbb{R}^{n},  \label{isotropic}
\end{equation}%
where $|X|$ denotes the euclidian norm of $X$.

A random complex vector $Y\in \mathbb{C}^{n}$ is called isotropic if $(\Re
Y,\Im Y)\in \mathbb{R}^{2n}$ is isotropic and $\Re Y$ and $\Im
Y$ are independent.
\end{definition}

The definition implies that an isotropic complex random vector satisfies
that $\mathbf{E}\{(Y,X)\}=0$ and $\mathbf{E}\{|(Y,X)|^{2}\}=n^{-1}|X|^{2}$
for any $X\in \mathbb{C}^{n}$, but the reverse is not true.

Recall also that a function $f:\mathbb{R}^{n}\rightarrow \mathbb{R}$ is
called log-concave if for any $\theta \in \lbrack 0,1]$ and any $%
X_{1},X_{2}\in \mathbb{R}^{n}$, then $f\big(\theta X_{1}+(1-\theta )X_{2}%
\big)\geq f(X_{1})^{\theta }f(X_{2})^{1-\theta }$.

\begin{definition}
[Log-concave measure] A measure $\mu $ on $\mathbb{R}^{n}$ (or $\mathbb{C}%
^{n}$) is log-concave if for any measurable subsets $A,B$ of \thinspace\ $%
\mathbb{R}^{n}$ (or $\mathbb{C}^{n}$) and any $\theta \in \lbrack 0,1]$,
\begin{equation*}
\mu (\theta A+(1-\theta )B)\geq \mu (A)^{\theta }\mu (B)^{(1-\theta )}
\end{equation*}%
whenever $\theta A+(1-\theta )B=\{\theta X_{1}+(1-\theta )X_{2}\,:\,X_{1}\in
A,\;X_{2}\in B\}$ is measurable.
\end{definition}

\begin{remark}
If $Y$ is a random vector with a log-concave distribution, then any affine
image of $Y$ has a log-concave distribution. If $Y_{1}$ and $Y_{2}$ are
independent random vector with log-concave distributions, then the pair $%
(Y_{1},Y_{2})$ has a log-concave distribution and $Y_{1}+Y_{2}$ has a
log-concave distribution as well (see \cite{{Da-Ko-Ha:69}, {Bo:74},{Pr:71}}%
).  The Brunn-Minkowski inequality provides examples of log-concave
measures, that are the uniform Lebesgue measure on compact convex
subsets of $\mathbb{R}^{n}$ as well as their marginals. More
generally Borell's theorem \cite{Bo:75} characterizes the
log-concave measures that are not supported by any hyperplane as the
absolutely continuous measures (with respect to Lebesgue measure)
with a log-concave density. Note that the distribution of an
isotropic vector is not supported by any hyperplane.
\end{remark}

The paper is organized as follows. In Section 2 we present necessary
spectral and probabilistic facts, including recent ones on the isotropic
random vectors. Section 3 contains the proof of our main result (Theorem \ref%
{t:MPP}). Our proof is different from that of \cite{Ma-Pa:67}, and
follows essentially the scheme outlined in \cite{Pa:99}.

\section{Necessary Spectral and Probabilistic Facts}

We will begin by recalling several facts on the resolvent of real symmetric
(hermitian) matrices. Here and below $|...|$ denotes the euclidian (or
hermitian) norm of vectors and matrices.

\begin{lemma}
\label{l:res} Let $A$ be a real symmetric (hermitian) matrix and
\begin{equation}
G_{A}(z)=(A-z)^{-1},\;\Im z\neq 0,  \label{res}
\end{equation}%
be its resolvent. We have:

\begin{itemize}
\item[(i)]
\begin{equation}
|G_{A}(z)|\leq |\Im z|^{-1},  \label{Gn}
\end{equation}

\item[(ii)] if for $Y\in \mathbb{R}^{n}(\mathbb{C}^{n})$ $L_{Y}$ is the
corresponding rank-one matrix defined in (\ref{LYX}), and $\tau \in \mathbb{R%
}$, then
\begin{equation}
G_{A+\tau L_{Y}}(z)=G_{A}(z)-\tau G_{A}(z)L_{Y}G_{A}(z)(1+\tau
(G_{A}(z)Y,Y))^{-1},\;\Im z\neq 0.  \label{rr1}
\end{equation}
\end{itemize}
\end{lemma}

The proof of Lemma \ref{l:res} is elementary. Next is a version of
the martingal-difference bounds for the variance of a function of
independent random variables (see \cite{ES:81}). The technique of
martingal differences was used in studies of random matrices by
Girko (see e.g. \cite{Gi:01}) for results and references).

\begin{lemma}
\label{l:mart} Let $\{Y_{\alpha }\}_{\alpha =1}^{m}$ be \ a collection of
i.i.d. random vectors of $\mathbb{R}^{n}(\mathbb{C}^{n})$ with a common
probability law $F$, and $\Phi :\mathbb{R}^{nm}(\mathbb{C}^{nm})\rightarrow
\mathbb{C}$ be a bounded borelian function, satisfying the inequalities:%
\begin{equation}
\sup_{X_{1},...,X_{m}\in \mathbb{R}^{n}(\mathbb{C}^{n})}|\Phi -\left. \Phi
\right\vert _{X_{\alpha }=0}|\leq C<\infty ,\;\alpha =1,...,m.  \label{FFa}
\end{equation}%
Then%
\begin{equation}
\mathbf{Var}\{\Phi (Y_{1},...,Y_{m})\}\leq 4C^{2}m  \label{boF}
\end{equation}
\end{lemma}

\begin{proof}
Denote for $\alpha =1,...,m-1$%
\begin{equation*}
\overline{\Phi }_{\alpha }:=\mathbf{E}\{\Phi |Y_{1},...,Y_{\alpha }\}=\int_{%
\mathbb{R}^{\alpha n}}\Phi (Y_{1},...,Y_{\alpha },Y_{\alpha
+1},...,Y_{m})F(dY_{\alpha +1})...F(dY_{m}),
\end{equation*}%
and $\Phi _{0}=\mathbf{E}\{\Phi \},\;\Phi _{m}=\Phi $, and $\Delta _{\alpha
}=\overline{\Phi }_{\alpha }-\overline{\Phi }_{\alpha -1},\;\alpha =1,..,m$.
Then%
\begin{equation*}
\Phi -\mathbf{E}\{\Phi \}=\sum_{\alpha =1}^{m}\Delta _{\alpha },
\end{equation*}%
hence%
\begin{eqnarray*}
\mathbf{Var}\{\Phi \} &:&=\mathbf{E}\{|\Phi -\mathbf{E}\{\Phi
\}|^{2}\}=\sum_{\alpha =1}^{m}\mathbf{E}\{|\Delta _{\alpha }|^{2}\} \\
&&+\sum_{1\leq \alpha <\beta \leq m}^{m}\mathbf{E}\{\Delta _{\alpha }%
\overline{\Delta _{\beta }}+\Delta _{\beta }\overline{\Delta _{\alpha }}\}.
\end{eqnarray*}%
It follows then from the definition of $\Delta _{\alpha }$ that all the term
of the second sum on the r.h.s. are zero and the terms of the first sum are
bounded by $4C^{2}$ in view of condition (\ref{FFa}). Hence, we obtain (\ref%
{boF}).
\end{proof}

\bigskip

We will also use two recent fundamental results about log-concave measures.
The first is a large deviation inequality due to G. Paouris \cite{Pao:06}:

\begin{theorem}
\label{thm: paouris} There exists a positive constant $C>0$ such that for
any integer $n\geq 1$ and any isotropic random vector $Y\in \mathbb{R}^{n}$
with a log-concave distribution we have
\begin{equation}
\mathbf{P}\{|Y|\geq Ct\}\leq \exp (-t\sqrt{n}).  \label{equ:paouris}
\end{equation}%
for every $t\geq 1$.
\end{theorem}

The second is an inequality on the concentration of the euclidian norm. It
is a reformulation of the main result from \cite{Kl:06}:

\begin{theorem}
\label{thm:klartag} There exist positive constants $C,\alpha $ with $\alpha
<1$, such that for any integer $n\geq 1$ and any isotropic random vector $%
Y\in \mathbb{R}^{n}$ with a log-concave distribution we have
\begin{equation}
\mathbf{Var}\{\,|Y|^{2}\}\leq C/n^{\alpha }.  \label{equ:klartag}
\end{equation}
\end{theorem}

Weaker forms of (\ref{equ:klartag}) with a power of $\ln n$ instead of
$n^{\alpha}$ were first proved in \cite{Kl2:06} and later in \cite{FGP:06}.

We can now state one main technical tool about log-concave measures that we
will need.

\begin{lemma}
\label{lem: variance bound} Let $A$ be a complex matrix, with
operator norm $\|A\|\leq 1$ and let $Y\in \mathbb{R}^{n}$ (or
$\mathbb{C}^{n}$) be an isotropic random vector with a log-concave
distribution. Then
\begin{equation}
\mathbf{Var}\{(AY,Y)\}\leq 4C/n^{\alpha }  \label{equ: variance
bound}
\end{equation}%
where $C<\infty $ and $\alpha \in (0,1)$ are the constants from (\ref%
{equ:klartag}).
\end{lemma}

\begin{proof}
Let us first consider the case $Y\in \mathbb{R}^{n}$. \ Writing $A=R+iI\,$,
where $R$ and $I$ are hermitian and denoting $\xi =\xi _{1}+i\xi
_{2}=(RY,Y)+i(IY,Y),\;\overset{\circ }{\xi }_{1,2}=\xi _{1,2}-\mathbf{E}%
\{\xi _{1,2}\}$, we have%
\begin{equation*}
\mathbf{Var}\{(AY,Y)\}:=\mathbf{E}\{\big|\overset{\circ }{\xi }\big|^{2}\}=\mathbf{E}%
\{(\overset{\circ }{\xi_{1}})^{2}\}+\mathbf{E}
\{(\overset{\circ }{\xi_{2}})^{2}\}=%
\mathbf{Var}\{{\xi }_{1}\}+\mathbf{Var}\{{%
\xi }_{2}\}.
\end{equation*}%
Hence the proof reduces to the case where $A$ is hermitian, for
which we have $(AY,Y)=(BY,Y)$ where $B=\Re A$ is a real symmetric matrix with $%
\|B\|\leq 1$. Thus we may assume without loss of \ generality that
$A$ is a real symmetric matrix. If $V$ is an isometry then $VY$ is
also an isotropic random vector with a log-concave distribution.
Hence, we can assume that $A$ is diagonal, i.e.,
$A=\mathrm{diag}\{a_{i}\}_{i=1}^{n}$ with $|a_{i}|\leq 1$ for all
$1\leq i\leq n$ and set
\begin{equation*}
\varphi (a_{1},\dots ,a_{n})=\mathbf{Var}\{(AY,Y)\}=\mathbf{E}\left\{
\left\vert \sum a_{i}y_{i}^{2}\right\vert ^{2}\right\} -\left\vert
n^{-1}\sum a_{i}\right\vert ^{2},
\end{equation*}%
where $Y=\{y_{i}\}_{i=1}^{n}$. Since $\varphi $ is a positive quadratic
form, its maximum on the cube $[-1,1]^{n}$ is attained on one of its
vertices $\{-1,1\}^{n}$. In order to estimate $\varphi $ on a vertex, let $I$
be a subset of $\{1,\dots ,n\}$. Then
\begin{equation*}
\mathbf{Var}\left\{ \sum_{i\in I}y_{i}^{2}-\sum_{i\notin I}y_{i}^{2}\right\}
\leq \left( \mathbf{Var}^{1/2}\left\{ \sum_{i\in I}y_{i}^{2}\right\} +%
\mathbf{Var}^{1/2}\left\{ \sum_{i\notin I}y_{i}^{2}\right\} \right) ^{2},
\end{equation*}%
and we get
\begin{equation*}
\mathbf{Var}\{(AY,Y)\}\leq 4\max_{I}\mathbf{Var}\left\{ \sum_{i\in
I}y_{i}^{2}\right\} ,
\end{equation*}%
where the maximum is taken over all subset $I\subset \{1,\dots ,n\}$.

Now observe that the projection $Y_{I}=(y_{i})_{i\in I}$ of $Y$ onto $%
\mathbb{R}^{I}$ is clearly an isotropic random vector and that by
\cite{Da-Ko-Ha:69} it has a log-concave distribution. Thus Klartag's
result \cite{Kl:06} (see also \cite{FGP:06}) can be applied and
gives that there exist $\alpha \in (0,1)$ such that
\begin{equation*}
\mathbf{Var}\{|Y_{I}|^{2}\}=\mathbf{E}\left\{ \sum_{I}y_{i}^{2}\right\}
-n^{-2}|I|^{2}\leq cn^{-2}|I|^{2-\alpha }\leq c/n^{\alpha },
\end{equation*}%
where $c,\alpha $ are universal constant (do no depend on $Y$, $n,|I|$) and $%
|I|$ denotes $\mathrm{Card}I$. We conclude that
\begin{equation*}
\mathbf{Var}\{(AY,Y)\}\leq 4c/n^{\alpha }.
\end{equation*}%
This complete the proof when $Y$ is real; the complex case is
similar.
\end{proof}

\section{Proof of the Main Result}

We will prove first auxiliary facts that can be of independent interest.

\begin{proposition}
\label{p:sa} Let $N_{n,m}$ be the Normalized Counting Measure of eigenvalues
of (\ref{Hnm}), in which $\{Y_{\alpha }\}_{\alpha =1}^{m}$ are i.i.d. random
vectors (not necessarily isotropic and/or with log-concave distribution) and $%
\{\tau _{\alpha }\}_{\alpha =1}^{m}$ are i.i.d. random variables,
independent of $\{Y_{\alpha }\}_{\alpha =1}^{m}$. Denote
\begin{equation}
g_{n,m}(z)=\int_{\mathbb{R}}\frac{N_{n,m}(d\lambda )}{\lambda -z},\;\Im
z\neq 0,  \label{gn}
\end{equation}%
the Stieltjes transform of $N_{n,m}$. Then we have:%
\begin{equation}
\mathbf{Var}\{N_{n,m}(\Delta )\}\leq 4m/n^2  \label{boN}
\end{equation}%
for any interval $\Delta \subset \mathbb{R}$, and
\begin{equation}
\mathbf{Var}\{g_{n,m}(z)\}\leq 4m/n^2|\Im z|^{2}  \label{bog}
\end{equation}%
for any $z\in \mathbb{C}\setminus \mathbb{R}$.
\end{proposition}

\begin{proof}
To prove (\ref{boN}) we use Lemma \ref{l:mart} with $\Phi =nN_{n,m}(\Delta )$%
, the number of eigenvalues of $H_{nm}$ in $\Delta $. Since%
\begin{equation*}
H_{n,m}-\left. H_{n,m}\right\vert _{Y_{\alpha }=0}=\tau _{\alpha
}L_{Y_{\alpha }}
\end{equation*}%
is a rank 1 matrix, we have by the mini-max principle%
\begin{equation*}
|nN_{n,m}-\left. nN_{n,m}\right\vert _{Y_{\alpha }=0}|\leq 1,
\end{equation*}%
i.e. the constant $C$ in (\ref{FFa}) is 1 in this case. This and (\ref{boF})
lead to (\ref{boN}).

In the case of $g_{n,m}$ we choose $\Phi =ng_{n,m}(z)$. Then we have
according to (\ref{NCM}), (\ref{gn}), and the spectral theorem for real
symmetric (hermitian) matrices%
\begin{equation*}
ng_{n,m}(z)=\mathrm{Tr}\,(H_{n,m}-z)^{-1}:=\mathrm{Tr}\,G_{n,m}(z).
\end{equation*}%
Then (\ref{rr1}) implies
\begin{equation}
|ng_{n,m}(z)-\left. ng_{n,m}(z)\right\vert _{Y_{\alpha }=0}|\leq \frac{|\tau
_{\alpha }||(G_{\alpha }^{2}Y_{\alpha },Y_{\alpha })|}{|1+\tau _{\alpha
}(G_{\alpha }Y_{\alpha },Y_{\alpha })|}.  \label{difg}
\end{equation}%
By spectral theorem for real symmetric (hermitian) matrices there exists
non-negative measure $m_{\alpha }$ such that for any integer $p$
\begin{equation*}
(G_{\alpha }^{p}Y_{\alpha },Y_{\alpha })=\int_{\mathbb{R}}\frac{m_{\alpha
}(d\lambda )}{(\lambda -z)^{p}}.
\end{equation*}%
Thus%
\begin{equation*}
|\tau _{\alpha }||(\pi G_{\alpha }^{2}Y_{\alpha },Y_{\alpha })|\leq |\tau
_{\alpha }|\int_{\mathbb{R}}\frac{m_{\alpha }(d\lambda )}{|\lambda -z|^{2}}%
,\;
\end{equation*}%
and%
\begin{equation*}
|1+\tau _{\alpha }(G_{\alpha }Y_{\alpha },Y_{\alpha })|\geq |\tau _{\alpha
}||\Im (G_{\alpha }Y_{\alpha },Y_{\alpha })|=|\tau _{\alpha }||\Im z|\int_{%
\mathbb{R}}\frac{m_{\alpha }(d\lambda )}{|\lambda -z|^{2}}.
\end{equation*}%
This implies
\begin{equation}
|ng_{n,m}(z)-\left. ng_{n,m}(z)\right\vert _{Y_{\alpha }=0}|\leq |\Im
z|^{-1}.  \label{dgga}
\end{equation}%
Thus we can choose $|\Im z|^{-1}$ as $C$ in (\ref{boF}) and obtain (\ref{bog}%
) from (\ref{boF}).
\end{proof}

\begin{proposition}
\label{p:MPE} Let $N^{(0)}$ and $\sigma $ be probability measures on $%
\mathbb{R}$ and $f^{(0)}$ is the Stieltjes transform (\ref{f0}) of $N^{(0)}$%
. Consider a non-negative measure $N$ on $\mathbb{R}$ and assume that its
Stieltjes transform (\ref{f}) satisfies (\ref{MPE}). We have

(i) if $\sigma $ in (\ref{MPE}) is such that
\begin{equation}
\overline{\tau }=\int_{\mathbb{R}}|\tau |\sigma (d\tau )\leq \infty ,
\label{expta}
\end{equation}%
then $N$ is a probability measure and is determined uniquely by (\ref{MPE});

(ii) if $N$ is a probability measure (i.e., $N(\mathbb{R)}=1$), then it is
uniquely determined by (\ref{MPE}) with any $\sigma $ (i.e. not necessarily
satisfying (\ref{expta})).
\end{proposition}

\begin{proof}
(i). Note that $N\neq 0$, otherwise $f=0$ and (\ref{MPE}) has the form%
\begin{equation}
0=f^{(0)}(z-c\overline{\tau }),\;\Im z\neq 0  \label{fest0}
\end{equation}%
which is impossible since $N^{(0)}\neq 0$ (see e.g. (\ref{SP})). Let us show
that $N$ is a probability measure, i.e., that $N(\mathbb{R)}=1$. It suffices
to show that
\begin{equation}
\lim_{y\rightarrow \infty }y|f(iy)|=1.  \label{lyf1}
\end{equation}%
(see \cite{Ak-Gl:93}, Section 59). Note that%
\begin{equation}
y\Im f(iy)=\int_{\mathbb{R}}\frac{y^{2}}{\lambda ^{2}+y^{2}}N(d\lambda
)\rightarrow N(\mathbb{R})>0,\;y\rightarrow \infty .  \label{yimf}
\end{equation}%
Denote%
\begin{equation}
\zeta (z)=-c\int_{\mathbb{R}}\frac{\tau \sigma (d\tau )}{1+\tau f(z)}.
\label{zetas}
\end{equation}%
and prove that%
\begin{equation}
\lim_{y\rightarrow \infty }y^{-1}\zeta (iy)=0.  \label{lyzet}
\end{equation}%
Indeed, write%
\begin{equation*}
\zeta (iy)=-c\int_{|\tau |<\sqrt{y}}\frac{\tau \sigma (d\tau )}{1+\tau f(iy)}%
-c\int_{|\tau |\geq \sqrt{y}}\frac{\tau \sigma (d\tau )}{1+\tau f(iy)}%
:=I_{1}+I_{2}.
\end{equation*}%
We have for $y\rightarrow \infty $ in view of the bound $|f(iy)|\leq 1/y$:%
\begin{equation*}
|I_{1}|\leq \frac{c\sqrt{y}}{1-\sqrt{y}/y}=O(\sqrt{y}),
\end{equation*}%
thus $\lim_{y\rightarrow \infty }y^{-1}I_{1}=0$. We have also%
\begin{equation*}
|I_{2}|\leq \frac{c}{\Im f(iy)}\int_{|\tau |\geq \sqrt{y}}\sigma (d\tau ),
\end{equation*}%
and in view of (\ref{yimf}) and the obvious relation%
\begin{equation*}
\lim_{y\rightarrow \infty }\int_{|\tau |\geq \sqrt{y}}\sigma (d\tau )=0
\end{equation*}%
we conclude that $\lim_{y\rightarrow \infty }y^{-1}I_{2}=0$. Thus, (\ref%
{lyzet}) is proved. This implies that $iy+\zeta
(iy)=iy(1+o(1)),\;y\rightarrow \infty $ and since $N^{(0)}$ is a probability
measure, we have
\begin{equation*}
\lim_{y\rightarrow \infty }y|f^{(0)}(iy)|=1.
\end{equation*}%
It follows then from (\ref{MPE}) that (\ref{lyf1}) is true, i.e., $N$ is a
probability measure.

Let us prove now that (\ref{MPE}) determines $N$ uniquely. Assume that there
exists two probability measures $N^{\prime }$ and $N^{\prime \prime }$,
whose Stieltjes transforms $f^{\prime }$ and $f^{\prime \prime }$ satisfy (%
\ref{MPE}). Since $f^{\prime }$ and $f^{\prime \prime }$ are analytic in the
upper halfplane, there exists a sequence $\{z_{j}\}_{j\geq 1}$ such that $%
z_{j}\rightarrow \infty $ as $j\rightarrow \infty $ and $(f^{\prime
}-f^{\prime \prime })(z_{j})\neq 0$. Subtracting (\ref{MPE}) for $f^{\prime
} $ from that for $f^{\prime \prime }$ we find%
\begin{equation*}
1=c\int_{\mathbb{R}}\frac{N^{(0)}(d\lambda )}{(\lambda -z-\zeta ^{\prime
}(z_{j}))(\lambda -z-\zeta ^{\prime \prime }(z_{j}))}\int_{\mathbb{R}}\frac{%
\tau ^{2}\sigma (d\tau )}{(1+\tau f^{\prime }(z_{j}))(1+\tau f^{\prime
\prime }(z_{j}))},
\end{equation*}%
where $\zeta ^{\prime }$ and $\zeta ^{\prime \prime }$ are given by (\ref%
{zetas}) for $f=f^{\prime },f^{\prime \prime }$. By using a bit more
sophisticated version of the proof of (\ref{lyf1}) it can be proved that the
limit of the r.h.s. of this relation is zero as $j\rightarrow \infty $. Thus
we have that $f^{\prime }=f^{\prime \prime }$. Since the Stieltjes transform
of a non-negative measure determines the measure uniquely (see e.g. (\ref{SP}%
)) we conclude that $N^{\prime }=N^{\prime \prime }$.

(ii). If $N$ is a probability measure, then we have $N(\mathbb{R})=1$ in the
r.h.s. of (\ref{yimf}). Now it is easy to check that the argument proving
assertion (i) is applicable to the case, where (\ref{expta}) is not valid
and we obtain the proof of assertion (ii) of the proposition.
\end{proof}

\bigskip

Now we are ready to prove

\begin{theorem}
\label{t:MPP} Let $\{Y_{\alpha }\}_{\alpha =1}^{m}$ be i.i.d. isotropic
random vectors of $\mathbb{R}^{n}($or $\mathbb{C}^{n})$ with a log-concave
distribution and $\{\tau _{\alpha }\}_{\alpha =1}^{m}$ be i.i.d. random
variables with a common probability law $\sigma $. Consider random matrices $%
H_{n,m}$ of (\ref{Mnm}) -- (\ref{Hnm}) and assume (\ref{N0N}). Then there
exist a probability measure $N$ such that for any interval $\Delta \subset
\mathbb{R}$ we have in probability%
\begin{equation}
\lim_{n\rightarrow \infty ,m\rightarrow \infty ,m/n\rightarrow c\in \lbrack
0,\infty )}N_{n,m}(\Delta )=N(\Delta )  \label{NnmNP}
\end{equation}%
for the Normalized Counting Measure $N_{n,m}$ of (\ref{NCM}) of eigenvalues
of $H_{n,m}$

The limiting non-random measure $N$ is uniquely determined by equation (\ref%
{MPE}) for its Stieltjes transform (\ref{f}).
\end{theorem}

\begin{proof}
In view of (\ref{boN}) it suffices to prove that the expectations%
\begin{equation}
\overline{N}_{n,m}=\mathbf{E}\{N_{n,m}\}  \label{Nbar}
\end{equation}%
of the Normalized Counting Measure (\ref{NCM}) of eigenvalues of $H_{nm}$
will converge weakly to the measure, whose Stieltjes transform solves of (%
\ref{MPE}). Recall now that weak convergence of probability measures to a
probability measure $N$ is equivalent to the convergence of their Stieltjes
transforms to the Stieltjes transform $f$ of (\ref{f}) $N$ on a compact set
of $\mathbb{C}\setminus \mathbb{R}$ and to the relation%
\begin{equation}
\lim_{y\rightarrow \infty }y|f(iy)|=1.  \label{mes1}
\end{equation}%
(see e.g. \cite{Ak-Gl:93}, Section 59). We are going to prove below these
two \ facts.

We derive first equation (\ref{MPE}). Assume first that the random variables
$\tau _{\alpha }$ are bounded:%
\begin{equation}
|\tau _{\alpha }|\leq T<\infty .  \label{tauT}
\end{equation}%
Write the resolvent identity for the resolvents%
\begin{equation}
G=(H_{n,m}-z)^{-1},\;\;\mathcal{G}=(H_{n}^{(0)}-z)^{-1}  \label{GG}
\end{equation}%
(we will omit below the sub-indices $n$ and $m$ and the argument $z$ in the
resolvents). We have in view of (\ref{Hnm}):%
\begin{equation*}
G=\mathcal{G-}\sum_{\alpha =1}^{m}\tau _{\alpha }GL_{Y_{\alpha }}\mathcal{G},
\end{equation*}%
and if
\begin{equation*}
\overline{G}=\mathbf{E}\{G\},
\end{equation*}%
then%
\begin{equation}
\overline{G}=\mathcal{G-}\sum_{\alpha =1}^{m}\mathbf{E}\{\tau _{\alpha
}GL_{Y_{\alpha }}\}\mathcal{G}.  \label{GG1}
\end{equation}%
Choose $t\geq 1$ and write (\ref{GG1}) as%
\begin{equation}
\overline{G}=\mathcal{G-}\sum_{\alpha =1}^{m}\mathbf{E}\{\tau _{\alpha
}GL_{Y_{\alpha }}\mathbf{1}_{|Y_{\alpha }|<t}\}\mathcal{G}-R_{1},
\label{GG2}
\end{equation}%
where%
\begin{equation}
R_{1}=\sum_{\alpha =1}^{m}\mathbf{E}\{\tau _{\alpha }GL_{Y_{\alpha }}\mathbf{%
1}_{|Y_{\alpha }|\geq t}\}\mathcal{G}.  \label{R1}
\end{equation}%
Denote%
\begin{equation}
G_{\alpha }=\left. G\right\vert _{Y_{\alpha }=0}.  \label{Gal}
\end{equation}%
Then we have from (\ref{rr1})%
\begin{equation*}
GL_{Y_{\alpha }}=G_{\alpha }L_{Y_{\alpha }}(1+\tau _{\alpha }(G_{\alpha
}Y_{\alpha },Y_{\alpha }))^{-1}.
\end{equation*}%
This allows us to write (\ref{GG2}) as%
\begin{equation*}
\overline{G}=\mathcal{G-}\sum_{\alpha =1}^{m}\mathbf{E}\left\{ \frac{\tau
_{\alpha }}{1+\tau _{\alpha }(G_{\alpha }Y_{\alpha },Y_{\alpha })}G_{\alpha
}L_{Y_{\alpha }}\mathbf{1}_{|Y_{\alpha }|<t}\right\} \mathcal{G}+R_{1}
\end{equation*}%
and then, denoting%
\begin{equation}
f_{n,m}(z):=\int_{\mathbb{R}}\frac{\overline{N}_{n,m}(d\lambda )}{\lambda -z}%
=\mathbf{E}\{n^{-1}\mathrm{Tr}(H_{n,m}-z)^{-1}\},\;\Im z\neq 0,  \label{fnm}
\end{equation}%
as%
\begin{equation}
\overline{G}=\mathcal{G-}\frac{m}{n}\int_{\mathbb{R}}\frac{\tau _{\alpha
}\sigma (d\tau )}{1+\tau f_{n,m}(z)}\overline{G}\mathcal{G}%
-\sum_{q=1}^{6}R_{q},  \label{GG3}
\end{equation}%
where $R_{1}$ is given by (\ref{R1}) and
\begin{equation}
R_{2}=\sum_{\alpha =1}^{m}\mathbf{E}\left\{ \left( \frac{\tau _{\alpha }}{%
1+\tau _{\alpha }(G_{\alpha }Y_{\alpha },Y_{\alpha })}-\frac{\tau _{\alpha }%
}{1+\tau _{\alpha }f_{n,m}}\right) G_{\alpha }L_{Y_{\alpha }}\mathbf{1}%
_{|Y_{\alpha }|<t}\right\} \mathcal{G},  \label{R2}
\end{equation}%
\begin{equation}
R_{3}=\sum_{\alpha =1}^{m}\mathbf{E}\left\{ \frac{\tau _{\alpha }}{1+\tau
_{\alpha }f_{n,m}}\left( G_{\alpha }L_{Y_{\alpha }}-n^{-1}G_{\alpha }\right)
\mathbf{1}_{|Y_{\alpha }|<t}\right\} \mathcal{G},  \label{R3}
\end{equation}%
\begin{equation}
R_{4}=\frac{1}{n}\sum_{\alpha =1}^{m}\mathbf{E}\left\{ \frac{\tau _{\alpha }%
}{1+\tau _{\alpha }f_{n,m}}\left( G_{\alpha }-G\right) \mathbf{1}%
_{|Y_{\alpha }|<t}\right\} \mathcal{G},  \label{R4}
\end{equation}%
\begin{equation}
R_{5}=\frac{m}{n}\mathbf{E}\left\{ \frac{1}{m}\sum_{\alpha =1}^{m}\left(
\frac{\tau _{\alpha }}{1+\tau _{\alpha }f_{n,m}}-\int_{\mathbb{R}}\frac{\tau
\sigma (d\tau )}{1+\tau f_{n,m}(z)}\right) G\mathbf{1}_{|Y_{\alpha
}|<t}\right\} \mathcal{G},  \label{R5}
\end{equation}%
\begin{equation}
R_{6}=\frac{m}{n}\int_{\mathbb{R}}\frac{\tau \sigma (d\tau )}{1+\tau
f_{n,m}(z)}\mathbf{E}\left\{ G\mathbf{1}_{|Y_{\alpha }|\geq t}\right\} .
\label{R6}
\end{equation}%
It follows from (\ref{GG3}) that if
\begin{equation}
\widetilde{\mathcal{G}}=\mathcal{G}(\widetilde{z}_{n,m}(z)),\;\widetilde{z}%
_{n,m}(z)=z-c\int_{\mathbb{R}}\frac{\tau \sigma (d\tau )}{1+\tau f_{n,m}(z)},
\label{GGzt}
\end{equation}%
then%
\begin{equation}
\overline{G}=\widetilde{\mathcal{G}}-\sum_{q=1}^{6}\widetilde{R}_{q},
\label{GG4}
\end{equation}%
where $\widetilde{R}_{q},\;q=1,...,6$ is obtained from $R_{q},\;q=1,...,6$
by replacing $\mathcal{G}$ by $\widetilde{\mathcal{G}}$ in (\ref{R1}), (\ref%
{R2}) -- (\ref{R6}). Note that
\begin{equation*}
\Im \widetilde{z}_{n,m}(z)=\Im z+c\Im f_{n,m}(z)\int_{\mathbb{R}}\frac{\tau
^{2}\sigma (d\tau )}{|1+\tau f_{n,m}(z)|},
\end{equation*}%
and since $\Im f_{n,m}(z)\Im z>0,\;\Im z\neq 0$ by (\ref{fnm}) we have%
\begin{equation*}
|\Im \widetilde{z}_{n,m}(z)|\geq |\Im z|
\end{equation*}%
and $\widetilde{\mathcal{G}}(z)$ is well defined for $\Im z\neq 0$.

Applying to (\ref{GG3}) the operation $n^{-1}\mathrm{Tr}$ and recalling (\ref%
{NCM0}), (\ref{GG}), (\ref{fnm}), and the spectral theorem, we obtain%
\begin{equation}
f_{n,m}(z)=f_{n}^{(0)}(\widetilde{z}_{n,m}(z))-\sum_{q=1}^{6}\widetilde{r}%
_{q},  \label{MPEnm}
\end{equation}%
where%
\begin{equation*}
f_{n}^{(0)}(z)=n^{-1}\mathrm{Tr}\mathcal{G}(z)=\int_{\mathbb{R}}\frac{%
N_{n}^{(0)}(d\lambda )}{\lambda -z},
\end{equation*}%
and
\begin{equation*}
\widetilde{r}_{q}(z)=n^{-1}\mathrm{Tr}\widetilde{R}_{q}.
\end{equation*}%
We will prove now that if
\begin{equation}
|\Im z|\geq 2t^{2}T,  \label{ztT}
\end{equation}%
where $t\geq 1$ and $T$ is defined in (\ref{tauT}), then there exist $%
C_{t,T}<\infty $ and $b>0$ such that%
\begin{equation}
|\widetilde{r}_{q}(z)|=o(1),\;q=1,...,6,  \label{rtq}
\end{equation}%
and we write here and below $o(1)$ for a quantity that tends to zero under
condition (\ref{nmc}) and (\ref{ztT}).

We begin from $R_{1}$. Note first that in view of (\ref{fnm}) $\Im
f_{n,m}(z)\Im z\geq 0,\;\Im z\neq 0$, hence%
\begin{equation}
|\widetilde{z}_{n,m}(z)|\geq |\Im z|\geq 2t^{2}T.  \label{ztnm}
\end{equation}%
This, (\ref{R1}), (\ref{ztT}), and (\ref{equ:paouris}) imply%
\begin{eqnarray*}
|\widetilde{r}_{1}(z)| &=&\left\vert \frac{1}{n}\sum_{\alpha =1}^{m}\mathbf{E%
}\{\tau _{\alpha }(\widetilde{\mathcal{G}}GY_{\alpha },Y_{\alpha })\mathbf{1}%
_{|Y_{\alpha }|\geq t}\}\right\vert  \\
&\leq &\frac{m}{4nt^{4}T^{2}}\mathbf{E}\{|Y_{\alpha }|^{2}\mathbf{1}%
_{|Y_{\alpha }|\geq t}\}=o(1).
\end{eqnarray*}%
Next, we have from (\ref{R2})
\begin{equation*}
\widetilde{r}_{2}(z)=\frac{1}{n}\sum_{\alpha =1}^{m}\mathbf{E}\left\{ \frac{%
(\tau _{\alpha })^{2}((G_{\alpha }Y_{\alpha },Y_{\alpha })-f_{n,m})}{(1+\tau
_{\alpha }(G_{\alpha }Y_{\alpha },Y_{\alpha }))(1+\tau _{\alpha }f_{n,m})}(%
\widetilde{\mathcal{G}}GY_{\alpha },Y_{\alpha })\mathbf{1}_{|Y_{\alpha
}|<t}\right\} .
\end{equation*}%
According to (\ref{Gn}) and (\ref{fnm}) we have $|(G_{\alpha }Y_{\alpha
},Y_{\alpha })|\leq |\Im z|^{-1}|Y_{\alpha }|^{2},\;|f_{n,m}|\leq |\Im
z|^{-1}$. Thus (\ref{ztT}) and (\ref{ztnm}) yield for $|Y_{\alpha }|<t$:
\begin{eqnarray*}
|1+\tau _{\alpha }(G_{\alpha }Y_{\alpha },Y_{\alpha })| &\geq &1-Tt^{2}/|\Im
z|\geq 1/2, \\
\;|1+\tau _{\alpha }f_{n,m}| &\geq &1-T/|\Im z|\geq 1/2.
\end{eqnarray*}%
This, (\ref{Gn}), and (\ref{ztnm}) lead to the bound%
\begin{equation}
|\widetilde{r}_{2}(z)|\leq \frac{4mt^{2}T^{2}}{n|\Im z|^{2}}\mathbf{E}%
\left\{ \mathbf{E}_{\alpha }\{|(G_{\alpha }Y_{\alpha },Y_{\alpha
})-f_{n,m}|\}\right\} ,  \label{rt22}
\end{equation}%
where $\mathbf{E}_{\alpha }\{...\}$ denotes the expectation only with
respect to $Y_{\alpha }$. Since $Y_{\alpha }$ is isotropic and $G_{\alpha }$
does not depend on $Y_{\alpha }$, we have (see (\ref{Gal}))%
\begin{equation*}
\mathbf{E}_{\alpha }\{(G_{\alpha }Y_{\alpha },Y_{\alpha })\}=n^{-1}\mathrm{Tr%
}G_{\alpha }:=g_{\alpha }.
\end{equation*}%
Besides, we have $\mathbf{E}\{g_{n,m}\}=f_{n,m}$ by definitions (\ref{NCM}),
(\ref{gn}), and (\ref{fnm}). Thus%
\begin{equation*}
\mathbf{E}\left\{ |\mathbf{E}_{\alpha }\{|(G_{\alpha }Y_{\alpha },Y_{\alpha
})-f_{n,m}|\}\right\} \leq \mathbf{E}\left\{ \mathbf{Var}_{\alpha
}^{1/2}\{(G_{\alpha }Y_{\alpha },Y_{\alpha })\}\right\} +\mathbf{E}\left\{
\mathbf{|}g_{\alpha }-g_{n,m}|\right\} +\mathbf{Var}_{\alpha
}^{1/2}\{g_{n,m}\}.
\end{equation*}%
Now (\ref{equ:klartag}) and (\ref{Gn}) yield the bound $2C^{1/2}/|\Im
z|n^{a/2}$ for the first \ term of the r.h.s., (\ref{dgga}) yields the bound
$1/n|\Im z|$ for the second term, and (\ref{bog}) yields the bound $%
2m^{1/2}/n|\Im z|$ for the third term. It follows then from
(\ref{rt22}) that $\widetilde{r}_{2}(z)=o(1)$.

Write now%
\begin{equation*}
\widetilde{r}_{3}(z)=\frac{1}{n}\sum_{\alpha =1}^{m}\mathbf{E}\left\{ \frac{%
\tau _{\alpha }}{1+\tau _{\alpha }f_{n,m}}\left( (\widetilde{\mathcal{G}}%
G_{\alpha }Y_{\alpha },Y_{\alpha })-n^{-1}\mathrm{Tr}\widetilde{\mathcal{G}}%
G_{\alpha }\right) \mathbf{1}_{|Y_{\alpha }|<t}\right\} ,
\end{equation*}%
and obtain similarly to the case of $\widetilde{r}_{2}$:%
\begin{equation*}
|\widetilde{r}_{3}(z)|\leq \frac{2Tm}{n}\mathbf{E}\left\{ \mathbf{Var}%
_{\alpha }^{1/2}\{\mathbf{(}\widetilde{\mathcal{G}}G_{\alpha }Y_{\alpha
},Y_{\alpha })\}\right\} \leq \frac{2Tmc^{1/2}}{n|\Im z|n^{a/2}}%
=O(1/n^{a/2}).
\end{equation*}%
In the case of
\begin{equation*}
\widetilde{r}_{4}(z)=\frac{1}{n}\sum_{\alpha =1}^{m}\mathbf{E}\left\{ \frac{%
\tau }{1+\tau f_{n,m}}n^{-1}\mathrm{Tr}\widetilde{\mathcal{G}}\left(
G_{\alpha }-G\right) \mathbf{1}_{|Y_{\alpha }|<t}\right\}
\end{equation*}%
we use again (\ref{rr1}) to write%
\begin{equation*}
\widetilde{r}_{4}(z)=\frac{1}{n^{2}}\sum_{\alpha =1}^{m}\mathbf{E}\left\{
\frac{\tau _{\alpha }}{1+\tau _{\alpha }f_{n,m}}\frac{\tau _{\alpha
}(G_{\alpha }\widetilde{\mathcal{G}}G_{\alpha }Y_{\alpha },Y_{\alpha })}{%
(1+\tau _{\alpha }(G_{\alpha }Y_{\alpha },Y_{\alpha }))}\mathbf{1}%
_{|Y_{\alpha }|<t}\right\} .
\end{equation*}%
Thus, similarly to the case of $\widetilde{r}_{2}(z)$ we have%
\begin{equation*}
|\widetilde{r}(z)_{4}|\leq \frac{2Tm}{n^{2}|\Im z|}=O(1/n).
\end{equation*}%
Denote%
\begin{equation*}
\xi _{\alpha }=\frac{\tau _{\alpha }}{1+\tau _{\alpha }f_{n,m}}-\int_{%
\mathbb{R}}\frac{\tau \sigma (d\tau )}{1+\tau f_{n,m}(z)},\;\alpha =1,...,m.
\end{equation*}%
Then it follows from (\ref{R5}) that%
\begin{equation*}
\widetilde{r}_{5}(z)=\frac{m}{n^{2}}\mathbf{E}\left\{ \frac{1}{m}%
\sum_{\alpha =1}^{m}\xi _{\alpha }\mathrm{Tr}\widetilde{\mathcal{G}}G\;%
\mathbf{1}_{|Y_{\alpha }|<t}\right\}
\end{equation*}%
Since $\left\{ \xi _{\alpha }\right\} _{\alpha =1}^{m}$ is a collection of
i.i.d. random variables, and $\mathbf{E}\left\{ \xi _{\alpha }\right\} =0$
we can write in view of (\ref{Gn}), (\ref{ztnm}), and (\ref{tauT})%
\begin{equation*}
|\widetilde{r}_{5}(z)|\leq \frac{m}{n|\Im z|^{2}}\mathbf{E}^{1/2}\left\{
\left( \frac{1}{m}\sum_{\alpha =1}^{m}\xi _{\alpha }\right) ^{2}\right\}
\leq \frac{m^{1/2}}{n|\Im z|^{2}}\mathbf{E}^{1/2}\left\{ \xi
_{1}^{2}\right\} =O(1/n^{1/2}).
\end{equation*}%
Finally we have
\begin{equation*}
\widetilde{r}_{6}(z)=\frac{m}{n}\int_{\mathbb{R}}\frac{\tau \sigma (d\tau )}{%
1+\tau f_{n,m}(z)}\mathbf{E}\left\{ n^{-1}\mathrm{Tr}\widetilde{\mathcal{G}}%
G\;\mathbf{1}_{|Y_{\alpha }|\geq t}\right\} ,
\end{equation*}%
hence%
\begin{equation*}
|\widetilde{r}_{6}|\leq \frac{2T}{|\Im z|^{2}}e^{-t\sqrt{n}/C}.
\end{equation*}%
Since the collection $\{\overline{N}_{nm}\}$ consists of probability
measures, there exists a non-negative measure $N$ which is a limiting point
of the collection in the sense (\ref{nmc}) with respect to the vague
convergence. The Stieltjes transforms of the corresponding subsequence of
the collection converge to the Stieltjes transform $f$ of $N$ uniformly on a
compact $K\subset \mathbb{C}\setminus \mathbb{R}$. Choosing $K\subset \{z\in
\mathbb{C}:|\Im z|\geq 2t^{2}T\}$ and using (\ref{N0N}) and (\ref{rtq}) we
can make the limit (\ref{MPEnm}) along the subsequence for $|\Im z|\geq
2t^{2}T$. The limiting equation, i.e. (\ref{MPE}) for $|\Im z|\geq 2t^{2}T$
implies that $\Im f(z)\neq 0,\;\Im z\neq 0$, otherwise we would have (\ref%
{fest0}) for some $z_{0},\;|\Im z_{0}|\geq 2t^{2}T$. Thus the both
parts
of obtained (\ref{MPE}) are analytic for $\Im z\neq 0$ and the equation (\ref%
{MPE}) is valid everywhere in $\mathbb{C}\setminus \mathbb{R}$. By
Proposition \ref{p:MPE} (i) $N$ is uniquely determined by the equation. This
proves the theorem under conditions (\ref{tauT}).

Consider now a general case and introduce the truncated random variables
\begin{equation}
\tau _{\alpha }^{T}=\left\{
\begin{array}{cc}
\tau _{\alpha }, & |\tau _{\alpha }|<T, \\
0, & |\tau _{\alpha }|\geq T,%
\end{array}%
\right.   \label{tautr}
\end{equation}%
for any $T>0$. Denote%
\begin{equation*}
H_{n,m}^{T}=H_{n}^{(0)}+\sum_{\alpha =1}^{m}\tau _{\alpha }^{T}L_{Y_{\alpha
}}.
\end{equation*}%
Then%
\begin{equation*}
\mathrm{rank}(H_{n,m}-H_{n,m}^{T})\leq \mathrm{Card}\{\alpha \in \lbrack
1,m]:|\tau _{\alpha }|\geq T\},
\end{equation*}%
and if $N_{n,m}^{T}$ is the Normalized Counting Measure of eigenvalues of $%
H_{n,m}^{T}$ and $\overline{N}_{n,m}^{T}$ is its expectation, then the
mini-max principle implies for any interval $\Delta \subset \mathbb{R}$:%
\begin{equation*}
|\overline{N}_{n,m}(\Delta )-\overline{N}_{n,m}^{T}(\Delta )|\leq \mathbf{P}%
\{|\tau _{1}|\geq T\}.
\end{equation*}%
Hence any limiting point $N$ in the sense of (\ref{nmc}) of the collection $%
\{\overline{N}_{n,m}\}$ with respect to the vague convergence satisfies the
same inequality:%
\begin{equation}
|N(\Delta )-N^{T}(\Delta )|\leq \mathbf{P}\{|\tau _{1}|\geq T\},  \label{NNT}
\end{equation}%
where we took into account that according to the first part of the theorem
the sequence $\{\overline{N}_{n,m}^{T}\}$ converges weakly to the
probability measure $\overline{N}^{T}$ for any $0<T<\infty $. Choosing here $%
\Delta =\mathbb{R}$, we obtain that $N$ is a probability measure, i.e. $N$
is limiting point of $\{\overline{N}_{n,m}\}$ in the sense (\ref{nmc}) with
respect to the weak convergence.

Besides it follows from (\ref{NNT}) that $N$ is a weak limiting points of
the family $\{N^{T}\}_{T>0}$.

According to the first part of the theorem and (\ref{tautr})\ the Stieltjes
transform $f^{T}$ of $N^{T}$ is a unique solution of the functional equation

\begin{equation}
f^{T}(z)=f_{0}\left( z-c\int_{|\tau |<T}\frac{\tau \sigma (d\tau )}{1+\tau
f^{T}(z)}\right) .  \label{MPET}
\end{equation}%
If $\{N^{T_{i}}\}_{i\geq 1}$ is a subsequence converging weakly to $N$, then
$\{f^{T_{i}}\}_{i\geq 1}$ converges to the Stieltjes transform $f$ of $N$
uniformly on a compact set $K\subset \mathbb{C}\setminus \mathbb{R}$. Since $%
N(\mathbb{R})=1$, there exists $C>0$, such that%
\begin{equation*}
\min_{z\in K}\Im f(z)=C>0.
\end{equation*}%
This and the uniform convergence of $\{f^{T_{i}}\}_{i\geq 1}$ to $f$ on $K$
implies that we have for all sufficiently big $T_{i}$:%
\begin{equation*}
\min_{z\in K}\Im f^{T_{i}}(z)=C/2>0.
\end{equation*}%
Thus $|\tau /(1+\tau f^{T_{i}}(z))|\leq (\Im f^{T_{i}}(z))^{-1}\leq
2/C<\infty ,\;z\in K$. This allows us to make the limit $T_{i}\rightarrow
\infty $ in (\ref{MPET}), i.e. to obtain (\ref{MPE}) for $f$. According to
Proposition \ref{p:MPE} (ii) the equation is uniquely soluble for any
probability measure $\sigma $, hence the limit point $f$ of the sequence $%
\{f^{T}\}_{T>0}$ is unique. Since $f$ determines uniquely $N$ (see (\ref{SP}%
)), the theorem is proved.
\end{proof}

\begin{remark}
Our result can be used for another random matrix, defined via the
vectors $\{Y_{\alpha}\}_{a=1}^{m}$ as%
\begin{equation}
\mathcal{G}_{n,m}=\{(Y_{\alpha },Y_{\beta })\}_{a,\beta =1}^{m}.
\label{Gram}
\end{equation}%
It  can be called the Gram matrix of the collection $\{Y_{\alpha
}\}_{a=1}^{m}$.

Denoting $Y_{\alpha }=\{y_{\alpha j}\}_{j=1}^{n}$ we can view
$Y=\{y_{\alpha
j}\}_{\alpha ,j=1}^{m,n}$ as~a $m\times n$ random matrix. Then we have for $%
\mathcal{G}_{n,m}$
\begin{equation*}
\mathcal{G}_{n,m}=YY^{\ast },
\end{equation*}%
and for $M_{n,m}$ of (\ref{Mnm}) with $\tau _{\alpha }=1,\;\alpha =1,...,m$%
\begin{equation*}
M_{n,m}=Y^{\ast }Y.
\end{equation*}%
Assume that the law of $Y_{\alpha }$ is continuous. Then for $n>m$
the vectors  $\{Y_{\alpha }\}_{a=1}^{m}$ are linearly independent
with probability 1, all the eigenvalues of $\mathcal{G}_{n,m}$ are
strictly
positive and coincide with non-zero eigenvalues of $M_{n,m}$. $%
M_{n,m}$ has $n-m$ zero eigenvalues, whose eigenvectors form a
basis \ of the complement of the span of $\{Y_{\alpha
}\}_{a=1}^{m}$\ in $\mathbb{R}^{n} $ (or $\mathbb{C}^{n}$). Denote
$\widetilde{N}_{n,m}$ the Normalized
Counting Measure of eigenvalues of $\mathcal{G}_{n,m}$. Then we have%
\begin{equation*}
\widetilde{N}_{n,m}=-\frac{n-m}{m}\delta _{0}+\frac{n}{m}N_{n,m}.
\end{equation*}%
Hence, if $N$ is the limit of $N_{n,m}$ in the sense of
(\ref{nmc}), then the limit $\widetilde{N}$ of
$\widetilde{N}_{n,m}$ in the same sense also exists and is equal
to $\widetilde{N}=-(c^{-1}-1)\delta _{0}+c^{-1}N$. Since in this
case $N=(1-c)\delta _{0}+N^{\ast }$, where $N^{\ast }(d\lambda
)=\rho ^{\ast }(\lambda )d\lambda $, where the support of $\rho
^{\ast }$ is $[a_{-},a_{+}],\;a_{\pm }=(1\pm \sqrt{c})^{2}$, and
\begin{equation}
\rho ^{\ast }(\lambda )=\frac{1}{4\pi c\lambda
}\sqrt{(a_{+}-\lambda )(\lambda -a_{-})},\;\lambda \in \lbrack
a_{-},a_{+}],  \label{rMP}
\end{equation}%
(see \cite{Ma-Pa:67}), we conclude that $\widetilde{N}$ is
absolutely continuous and has the density $c^{-1}\rho ^{\ast }$.

Similar argument shows that for $m>n$ $\ $we have $\widetilde{N}%
=(1-c^{-1})\delta _{0}+c^{-1}N$, where now $N$ is absolutely
continuous and its density is (\ref{rMP}).
\end{remark}


\begin{thebibliography}{99}
\bibitem{Ak-Gl:93} N. I. Akhiezer,~I. M. Glazman, \emph{Theory of Linear
Operators in Hilbert Space}, Dover, New York, 1993

\bibitem{Au:06} G. Aubrun, Random points in the unit ball of $\ell_p^n$,
\emph{Positivity} \textbf{10} (2006) 755-759.

\bibitem{Au:07} G. Aubrun, Sampling convex bodies: a random matrix approach,
\emph{Proceedings AMS}, \textbf{135} (2007) 1293-1303.

\bibitem{Bo:74} C. Borell, Convex measures on locally convex spaces  \emph{%
Ark. Math.}, \textbf{12} (1974) 239-252.

\bibitem{Bo:75} C. Borell, Convex set functions in $d$-space \emph{Period.
Math. Hungar. }, \textbf{6} no.2, (1975) 111-136.

\bibitem{Bou:96} J. Bourgain, Random points in isotropic convex sets, In
\emph{Convex Geometric analysis}, (Berkeley, CA, 1996), vol. 34 of \emph{%
Math. Sci. Res. Inst. Publ.} Cambridge Univ. Press, Cambridge, 1999, pp.
53-58.

\bibitem{Da-Ko-Ha:69} Ju. S. Davidovi$\check{c}$, B. I. Korenbljum, B. I.
Hacet, A certain property of logarithmically concave functions. (Russian)
\emph{Dokl. Akad.Nauk SSSR}, \textbf{185} (1969) 1215-1218. \emph{Soviet
Math. Dokl. } \textbf{10} (1969) 447-480.

\bibitem{Da-Sz:01} K. R. Davidson, S. Szarek, Local operator theory, random
matrices and Banach spaces. In Handbook on the Geometry of Banach spaces,
Volume 1, W. B. Johnson, J. Lindenstrauss eds., Elsevier Science 2001, pp.
317-366.

\bibitem{ES:81} B. Efron, C. Stein. The jackknife estimate of variance.
\emph{Annals of statistics} \textbf{9} (1981) 586-596.

\bibitem{FGP:06} B. Fleury, O. Gu\'edon and G. Paouris,
A stability result for mean width of $L_p$-centroid bodies,
to appear in
Advances in Mathematics.

\bibitem{Gi-Mi:00} A. Giannopoulos, V. Milman, Concentration property on
probability spaces \emph{Advances in Mathematics} \textbf{156} (2000) 77-106.

\bibitem{Gi-Ts:03} A. Giannopoulos, A. Tsolomitis, Volume radius of a random
polytope in a convex body \emph{Math. Proc. Cambridge Philos. Soc.} \textbf{%
134} (2003) 13-21.

\bibitem{Gi:01} V. L. Girko, \emph{Theory of Stochastic Canonical Equations}%
, vols.I, II Kluwer, Dordrecht, 2001.


\bibitem{Kl:06} B. Klartag, Power-law estimates for the central limit
theorem for convex sets. arXiv:math.MG/0611577

\bibitem{Kl2:06} B. Klartag,
A central limit theorem for convex sets, Invent. Math., Vol. 168,
(2007), 91--131.

\bibitem{Li-Ru-Paj-Tom:05} A. Litvak, A. Pajor M. Rudelson, N.
Tomczak-Jaegermann, Smallest singular value of random matrices and geometry
of random polytopes \emph{Advances in Mathematics} \textbf{195} (2005)
491--523.

\bibitem{Ma-Pa:67} V. Marchenko, L. Pastur, The eigenvalue distribution in
some ensembles of random matrices, Math. USSR Sbornik \textbf{1} (1967)
457-483.

\bibitem{Mui:82} R. J. Muirhead, \emph{Aspects of Multivariate Statistical
Theory}, Wiley, New York, 1982.

\bibitem{Pao:06} G. Paouris, Concentration of mass of isotropic convex
bodies. C.R. Math. Acad. Sci. \textbf{342} (2006) 179--182

\bibitem{Pa:99} L. Pastur, A simple approach to the global regime of the
random matrix theory, In: \emph{Mathematical Results in Statistical Mechanics%
}, S. Miracle-Sole, J. Ruiz, V. Zagrebnov (Eds.), World Scientific,
Singapore 1999, pp. 429-454.

\bibitem{Pr:71} A. Pr\'ekopa, Logarithmic concave measures with application
to stochastic programming \emph{Acta Sci. Math. (Szeged)}, \textbf{32}
(1971) 301-316.

\bibitem{Sz:06} S. Szarek, Convexity, Complexity, and High Dimensions. In
"Proceedings of the International Congress of Mathematicians (Madrid,
2006)," Vol. II. European Math. Soc. 2006, 1599--1621.
\end{thebibliography}
\end{document}